\input amstex
\documentstyle{amsppt}

\def\C{\Gamma}
\def\A{\Cal A}
\def\B{\Cal B}

\magnification=1200

\topmatter
\title  {Group splittings and asymptotic topology }
\endtitle
\rightheadtext{ Group splittings and asymptotic topology }

\author{ Panos Papasoglu }
\endauthor

\date { }
\enddate
\affil {Universit\'e de Paris-Sud }\endaffil

\address{ 
Departement de Math\'ematiques, Universit\'e de Paris XI (Paris-Sud)
91405 Orsay, FRANCE.} 
\endaddress

\email{ panos\@topodyn.math.u-psud.fr } 
\endemail

\endsubjclass
\abstract { It is a consequence of the theorem of Stallings on groups with
many ends that splittings over finite groups are preserved by quasi-isometries.
In this paper
we use asymptotic topology to show that group splittings
are preserved by quasi-isometries in many cases. Roughly speaking we show that
splittings are preserved under quasi-isometries when the vertex groups are 
 fundamental groups of aspherical manifolds (or more generally `coarse $PD(n)$-groups')
and the edge groups are `smaller' than the vertex groups. }

\endabstract

\endtopmatter

\document

\NoBlackBoxes
\heading {\bf \S 0. Introduction } \endheading

The notion of quasi-isometry and the study  of the relation of
large-scale geometry of groups and algebraic properties has become
predominant in group theory after the seminal papers of Gromov \cite {G1,G2 }.\newline
A classical theorem of Stallings (\cite {St }) implies, 
that if $G$ splits over a finite group and $H$ is a group quasi-isometric
to $G$ then $H$ also splits over a finite group. \newline
Bowditch has shown recently (\cite {Bo}) that splittings of hyperbolic groups over 2-ended
groups are preserved by quasi-isometries.\newline
By a theorem of Kapovich and Leeb \cite {K-L} it follows
 that a group quasi-isometric to a non-geometric
Haken 3-manifold splits over a group commensurable to a surface group.
In this paper we show that group splittings are preserved by quasi-isometries
in many cases. 
Our approach  is based on asymptotic topology (`coarse topology') methods.
Schwartz's asymptotic version of the Jordan separation theorem (\cite {Sch,F-S}) is our main
tool. In fact we will need a stronger version of Schwartz's theorem that has been given recently
by Kapovich and Kleiner (\cite {K-K})
 in the context of their work on coarse $PD(n)$-spaces and groups. 
We will therefore
formulate our resluts in this more general setting.\newline
We say that a group $G$ is a coarse $PD(n)$ group if it acts discretely co-compactly
on a coarse $PD(n)$-space. We say that $G$ is a coarse $PD(n)$ group of dimension $n$
if it is a coarse $PD(n)$-group that has an $n$-dimensional $K(G,1)$.\newline
We note that examples
of coarse $PD(n)$-groups are fundamental groups of closed aspherical
$n$-manifolds. \newline
Using the theory of  coarse $PD(n)$-spaces we show that
 groups that are quasi-isometric to `trees of coarse $PD(n)$-spaces' split.
To pass from
geometry to algebra we use a recent result of Scott and Swarup (\cite {S-S}).\newline
We explain briefly the notation we use for graphs of groups. For more
details see \cite {Ba}, \cite {Se}. A graph of groups is given by
the following data:\newline
a) A finite graph $\C $. Each edge $e\in \C$ is oriented. We denote
by $\partial _0 e$ the initial vertex of $e$ and by $\partial _1 e$ 
the terminal vertex of $e$.\newline
b)To each vertex $v\in \C$ and each edge $e\in \C$ there correspond
groups $A_v,A_e$. If $v=\partial _0e$ or $v=\partial _1e$ we
have monomorphisms (respectively) $i_0:A_e\to A_v$ and $i_1:A_e\to A_v$.
We denote this collection of groups and morphisms by $\Cal A$.\newline 
Using these data one defines the fundamental group $\pi _1(\C,\A)$.
For the definition see \cite {Se}.\newline
With this notation we have the following:
\proclaim {Theorem 3.1 }
Let $G$ be a finitely generated group admitting a graph of groups decomposition
$G=\pi _1(\Gamma ,\A )$ such that all edge and vertex groups are coarse $PD(n)$ groups
of dimension $n$.
Suppose further that $(\Gamma,\A) $ is not a loop with all edge to vertex
maps isomorphisms and that it is not a graph of one edge with both edge to
vertex maps having as image an index 2 subgroup of the vertex group.
 If $H$ is quasi-isometric to $G$ 
then $H$ splits over a group that is quasi-isometric to an edge group of $(\Gamma,\A) $.
\endproclaim
We obtain the following corollaries:
\proclaim {Corollary 3.2  (\cite {F-M1,2}) }
Let $G$ be a solvable Baumslag-Solitar group.  If $H$
is a group quasi-isometric to $G$ then $H$ is commensurable
to a solvable Baumslag-Solitar group.
\endproclaim
\proclaim {Corollary 3.3 }
Let $G$ be a finitely generated group admitting a graph of groups decomposition
$G=\pi _1(\Gamma ,\A )$ such that all edge and vertex groups are virtually $\Bbb Z ^n$.
Suppose further that $(\Gamma,\A) $ is not a loop with all edge to vertex
maps isomorphisms and that it is not a graph of one edge with both edge to
vertex maps having as image an index 2 subgroup of the vertex group.
 If $H$ is quasi-isometric to $G$ 
then $H=\pi _1(\Delta ,\Cal B)$ where all vertex and edge groups of $\Delta $
are virtually $\Bbb Z ^n$.
\endproclaim
It turns out that splittings are invariant under quasi-isometries 
in the case that edge groups are `smaller' than vertex groups:
\proclaim {Theorem 3.4 }
Let $G$ be a finitely generated group admitting a graph of groups decomposition
$G=\pi _1(\Gamma ,\A )$ such that all vertex groups are coarse $PD(n)$-groups
 and all edge groups are dominated by coarse $PD(n-1)$ spaces.
 If $H$ is quasi-isometric to $G$ 
then $H$ splits over some group quasi-isometric to an edge group
of $\Gamma $.
\endproclaim
We say that a group $G$ is dominated by a coarse $PD(n)$-space $X$ if there is
a uniform embedding $f:G\to X$ (for more details see sec. 3). Note that a subgroup
of a coarse-$PD(n)$ group is dominated by a coarse $PD(n)$-space. So for example a free group
is dominated by a coarse $PD(2)$-space. \newline
The main geometric observation on which our results are based is that
the groups in theorems 3.1, 3.4 are `trees of spaces'. The simplest example of such
groups are products $\Bbb F _k\times \Bbb Z ^n$. To make the exposition
easily accessible to readers not familiar with the geometry of graphs
of groups and with `coarse $PD(n)$-spaces' we treat this special case first
in section 2. All `asymptotic topology' arguments that we need are already present
in this case. The link between algebra and geometry is provided by
a result of Scott and Swarup (\cite {S-S}) generalizing
the algebraic torus theorem of Dunwoody-Swenson (\cite { D-S }).\newline
 In section 3 we explain how to generalize these arguments 
to graphs of groups in which all edge and vertex groups
 are `coarse $PD(n)$-groups of dimension $n$'. For this it suffices
to understand the `tree-like shape' of graphs of
groups in general. The geometries of such groups have been described in several
places (see \cite {S-W}, \cite {Ep}, \cite {F-M1,2 }, \cite {Wh}).
We use similar arguments to treat the case of graphs of groups with vertex
groups coarse $PD(n)$ groups and edge groups, groups that are `dominated' by coarse
 $PD(n-1)$ groups.\newline
In section 4 we discuss how the results of this paper (and Stallings' theorem) could
be generalized and we ask some specific questions.\newline
In the course of this work we found out that some of our results had been
obtained earlier,independently, by Mosher, Sageev and Whyte. In particular they have
shown a stronger version of theorem 3.1, Corollary 3.4 and some cases of theorem 3.4 (\cite {MSW},\cite {MSW1}).
The main novelty (apart from the difference in the proofs)
of this paper compared to \cite {MSW} is that we improve $n-2$ to $n-1$ in theorem 3.4. So for
example from our results it follows that if $G,H$ are quasi-isometric groups and $G$ is an
amalgam of two aspherical 3-manifold
groups along a surface group  (or free group) then $H$ also splits over a virtual surface
(or virtually free) group. The work of \cite {MSW} implies a similar result when $G$ is
an amalgam over $\Bbb Z$.\newline
I would like to thank G.A. Swarup for his criticism of an earlier version of this paper
and B.Kleiner for explaining to me \cite {K-K}.

\heading {\bf \S 1. Preliminaries } \endheading 
A (K,L)-{\it quasi-isometry} between two metric spaces $X,Y$ is a map
$f:X\to Y$ such that the following two properties are satisfied:\newline
1) ${1\over K}d(x,y)-L\leq d(f(x),f(y))\leq Kd(x,y)+L$ for all $x,y\in X$.\newline
2) For every $y\in Y$ there is an $x\in X$ such that $d(y,f(x))\leq K$.\newline
We will usually simply say quasi-isometry instead of (K,L)-quasi-isometry.
Two metric spaces $X,Y$ are called quasi-isometric if there is a quasi-isometry
$f:X\to Y$.\newline
A {\it geodesic metric space }  is a metric space in which any
two points $x,y$ are joined by a path of length $d(x,y)$. In what
follows we will be interested in graphs that we always turn into
geodesic metric spaces by giving each edge length 1.\newline
We recall now some notation and results from \cite {K-K}. We refer the reader
to this paper for more details.\newline
Let $X$ be a connected locally finite simplicial complex. The 1-skeleton $X^1$
is a graph and we turn it into a geodesic metric space  as explained above.
If $A\subset X^1$ the $r$-neighborhood of $A$ is the set of points
in $X^1$ at distance less or equal to $r$ from $A$.
More generally if $K$ is a subcomplex of $X$ we define the $r$-neighborhood of $K$ ,$N_r(K)$,
to be the set of simplices intersecting $N_r(K^1)$.\newline
The {\it diameter }, $diam(K)$, of $K$ is by definition the diameter of $K^1$.\newline
We call a map $f:X\to Y$ between metric spaces a {\it uniform embedding} (see \cite {G2})
if the following two conditions are satisfied:\newline
1) There are $K,L$ such that for all $x,y\in X$ we have $d(f(x),f(y))<K(d(x,y))+L$.\newline
2) For every $E>0$ there is $D>0$ such that $diam (A)<E\Rightarrow diam(f^{-1}A)<D$.\newline
We say that the {\it distortion } of $f$ is bounded by $h$ ,where $h:\Bbb R ^+\to \Bbb R ^+$,
if for all $A\subset Y$, $diam(f^{-1}A)\leq h(diam (A))$.\newline
If $f$ satisfies only condition 1 above we say that $f$ is $(K,L)$-{\it lipschitz } map.\newline
$X$ is called {\it uniformly acyclic} if for each $R_1$ there is an $R_2$ such that
for each subcomplex $K\subset X$ with $diam(K)<R_1$ the inclusion $K\to N_{R_2}(K)$
induces zero on reduced homology groups.\newline
If $K\subset X$ is a subcomplex of $X$ and $R>0$ we say that a component of $X-N_R(K)$
is {\it deep } if it is not contained in $N_{R_1}(K)$ for any $R_1>0$.\newline
We say that $K$ {\it coarsely separates} $X$ if there is an $R>0$ such that
$X-N_R(K)$ has at least two deep components. \newline 
The appropriate context for
the results in this paper seems to be that of `coarse $PD(n)$ spaces'.
We refer to
 \cite {K-K} for a definition and an exposition of the theory of `coarse $PD(n)$ spaces'.
Important examples of `coarse $PD(n)$ spaces' are 
uniformly acyclic 
PL-manifolds of bounded geometry.\newline
We say that $X$ is a `coarse $n$-dimensional PL-manifold'  if
$X$ is quasi-isometric to a uniformly acyclic $n$-dimensional PL-manifold of bounded geometry. \newline
A reader not familiar with \cite {K-K} can read  this paper by replacing
everywhere `coarse $PD(n)$ space' by `coarse $n$-dimensional PL-manifold'. The only result
that we will use from \cite {K-K} is the coarse Jordan theorem stated below
(see theorem 7.7, footnote 11 and corollary 7.8 of \cite {K-K}). 
\proclaim {Proposition 1 (Coarse Jordan theorem) }
Let $X,X'$ be coarse $PD(n)$, $PD(n-1)$ spaces
respectively, $Z\subset X'$ and let $g:Z\to X$ be a uniform embedding such
that $g$ is a $(K,L)$ lipschitz map that is a uniform embedding with distortion
bounded by $f$.
Then \newline
1. If $Z=X'$ then there is an $R>0$ such that $X-N_R(g(X'))$ has exactly
2 deep components.\newline
2. There is an $N>0$ such that every non-deep component of $X-g(Z)$ is contained
in the $N$-neighborhood of $g(Z)$.\newline
3. There is an $M>0$ such that if $X-g(Z)$ has more than
one deep component, $X'$ is contained in the $M$-neighborhood of $Z$.\newline
4. If $Z=X'$ then for every $r$, each point of $N_r(g(X')) $ lies within uniform distance
from each of the deep components of $X-N_r(g(X'))$.\newline
The constants $M,N,R$ depend only on $K,L,f,X,Y$.
\endproclaim
We note that the metric on $Z$ in the above proposition is the metric induced
by $X'$.\newline
We say that $A\subset X$ {\it coarsely contains } $B$ if for some $R$, $B\subset N_R(A)$.
Let $f: X\to Y$ be a uniform embedding of a coarse $PD(n-1)$ space $X$ to a coarse
$PD(n)$ space $Y$. Let $K>0$ be such that $Y-N_K(f(X))$ has two deep components.
We will call a deep component of $Y-N_K(f(X))$ a {\it half coarse PD(n)-space }.\newline
We call a group $G$ a {\it coarse PD(n)-group} if it acts discretely co-compactly
simplicially on a coarse $PD(n)$ space.
We say that $G$ is a coarse $PD(n)$ group of { \it dimension } $n$
if it is a coarse $PD(n)$-group that has an $n$-dimensional $K(G,1)$.\newline
We call $(G;\{F_1,...,F_n\})$ a {\it coarse PD(n)-pair } if :\newline
1. $G$ acts discretely simplicially on a coarse $PD(n)$ space $X$ and there
is a connected $G$-invariant subcomplex $K\subset X$ such that the stabilizer
of each component of $X-K$ is conjugate to one of the $F_i's$.\newline
2. The $F_i$'s are coarse $PD(n-1)$ groups.

\heading {\bf \S 2. The geometry of direct products of abelian and free groups} \endheading

\proclaim {Theorem 2.1 }
Let $G=\Bbb F_k\times \Bbb Z ^n$ where $\Bbb F _k $ is the free group
on $k>1$ generators. If $H$ is quasi-isometric to $G$ 
then $H=\pi _1(\Delta ,\B )$ where all vertex and edge
groups of $\Delta $ are virtually $\Bbb Z ^n$.
\endproclaim
\demo{Proof}
Let $\Gamma _G$ be the Cayley graph of $G$ with respect to the standard
generators.  $\Gamma $ is isometric to $T_{2k}\times \Bbb R ^n$
where $T_{2k} $ is the homogeneous tree of degree $2k$. Clearly
$G$ acts discretely and co-compactly on $T_{2k}\times \Bbb R ^n$.

Let's denote $T_{2k}\times \Bbb R ^n$ by $X$ and let $p:X\to T_{2k}$ be the natural
projection from $X=T_{2k}\times \Bbb R ^n$ to $T_{2k}$.
 With this notation we have:
\proclaim {Lemma 2.2 }
Let $f:X\to X$ be a quasi-isometry. Then for any vertex $v$ of $T_{2k}$
there is a vertex $u\in T_{2k}$ such that $f(v\times \Bbb R ^n)$
and $u\times \Bbb R ^n $ are at finite distance from each other.

\endproclaim
\demo {Proof }
 
We note  that if $v$ is a vertex of $T_{2k}$ $v\times \Bbb R ^n$ 
separates
$T_{2k}\times \Bbb R ^n$ in more than 2 deep components. \newline
Moreover there are geodesic rays $r_1,r_2,r_3:[0,\infty )\to X$
such that $r_1,r_2,r_3$ lie in distinct components of $X-p^{-1}(v)$
and $d(r_1(t),p^{-1}(v))=t$.\newline
Let $K>0$ be such that $N_K(f(p^{-1}(v)))$ separates $X$
in more than 2 deep components. Clearly for $K$ sufficiently big $f(r_i),\, i=1,2,3$ are
 coarsely contained in distinct deep components
of $X-N_K(f(p^{-1}(v)))$. We pick $K$ so that this holds.
Let's call $C_i$ the deep component
coarsely containing $f(r_i)$.\newline
To simplify notation we set $S=N_K(f(p^{-1}(v)))$.\newline
$f(r_i)$ is not necessarily connected. Let $R_i$ be the path obtained
by joining $f(r_i(n))$ to $f(r_i(n+1))$ by a geodesic path for all $n\in \Bbb N$.
Clearly $R_i$ is coarsely contained in $C_i$.
 If we parametrize $R_i$ by arclength we
have that $d(R_i(t),S)$ is a proper function from $[0,\infty )$ to $[0,\infty )$.
Let $l_1=p(R_1)$, $l_2=p(R_2)$. We pick now a geodesic $l\subset T_{2k}$ such
that $p^{-1}l\cap R_1$ and $p^{-1}l\cap R_2$ are both unbounded. We explain
how to find such an $l$: If $l_1,l_2$ are finite then there are closed edges 
$e_1,e_2 $
of $T_{2k}$ such that $p^{-1}e_1\cap R_1, p^{-1}e_2\cap R_2 $ are both unbounded, so
we simply pick $l$ to be any geodesic containing both these edges. If both  $l_1,l_2$ are
infinite, since they are connected they contain at least one geodesic ray each.
We pick therefore $l$ to be a line  having an infinite intersection with
both geodesic rays. If one of them is infinite and one finite we pick $l$ similarly
, requiring that it has an infinite intersection with a given ray and passes
from a given edge. \newline
We have therefore that $S$ separates $p^{-1}l$. By part 3 of proposition 1 we
conclude that there is a $K_1$ such that $S$ is contained
in the $K_1$ neighborhood of $p^{-1}l$. \newline
We note that by parts 1,2 of proposition 1, $p^{-1}l\cap R_3$ is bounded.
We can therefore find a geodesic ray $r\in T_{2k}$ intersecting $l$
only at one point $u$, such
that $p^{-1}r\cap R_3$ is unbounded. We have then $l=r_1'\cup r_2'$ and $r_1'\cap r_2'=u$
where $r_1',r_2'$ are geodesic rays. As before we have that $S$ separates both
$p^{-1}(r\cup r_1')$ and $p^{-1}(r\cup r_2')$ so it lies in a finite neighborhood
of both. We conclude that $S$ lies in a finite neighborhood of all three:
$p^{-1}(r\cup r_1')$, $p^{-1}(r\cup r_2')$ and $p^{-1}l$. Therefore $S$ lies
in a finite neighborhood of $p^{-1}u$.\newline
We remark that if $f$ is an $(A,B)$-quasi-isometry then the proof above shows
that there is a $C$ that depends only on $A,B$ (and $X$) such that $f(v\times \Bbb R ^n)$
and $u\times \Bbb R ^n $ are in the $C$-neighborhood of each other. This is so
because the constants in proposition 1 depend only on $A,B$. 
We will use this fact in the next lemma.
$\blacksquare $
\enddemo
\proclaim {Lemma 2.3 }
Let $G=\Bbb F_k\times \Bbb Z ^n$ where $\Bbb F _k $ is the free group
on $k>1$ generators. If $H$ is quasi-isometric to $G$ 
then $H$ splits over a group commensurable to $\Bbb Z ^n$.
\endproclaim
\demo{Proof}

Let $\Gamma _H$ be the Cayley graph of $H$. If $v\in T_{2k}$ is a vertex
and $f:X\to \Gamma _H$ a quasi-isometry we have that $f(v\times \Bbb R ^n)$
coarsely separates $\Gamma _H$ to more than 2 deep components.\newline
Let $K$ be such that $\Gamma _H-N_K(f(v\times \Bbb R ^n))$ has more
than 2 deep components. We set $S=N_K(f(v\times \Bbb R ^n))$. Note
that $\Gamma _H-S$ has a finite number of deep components.\newline
We pick $M>0$ such that the following holds:\newline
For all $h\in H $ if $hS\cap S \ne \emptyset $ then $hS\subset N_M(S) $.\newline
It follows from our remark at the end of the proof of lemma 2.2 that such an $M$ exists.\newline
We will show that there is a subgroup $J$ of $H$ quasi-isometric to $S$
 such that $\Gamma _H/J$ has more than two ends.\newline
We fix a vertex $e\in S$. We define an equivalence relation on the
set of vertices $x\in S$:\newline
Let $g_x\in H$ such that $g_x x=e$. Let $C_1,...,C_k$ be the deep components of $\Gamma _H-S$
and let $D_1,...,D_m$ be the deep components of $\Gamma _H-N_M(S)$.\newline
We note that each $D_i$ ,$i=1,...,m$ is contained in some $C_j,\, j=1,...,k$.\newline
For each $g_x $ and $C_j$ we have that there are $i_1,...,i_r$ such that
$g_xC_j$ and $D_{i_1}\cup ...\cup D_{i_r}$ are at finite distance from each other.
We use this to define a map $f_x:\{C_1,...,C_k\}\to \Cal P (\{ D_1,...,D_m \})$ 
where $f_x(C_j)=\{D_{i_1}, ..., D_{i_r} \}$ if and only if
$g_xC_j$ and $D_{i_1}\cup ...\cup D_{i_r}$ are at finite distance from each other.\newline
Let's write now $x\sim y$ for $x,y$ vertices of $S$ if $f_x=f_y$. Clearly $\sim $
is an equivalence relation with finitely many equivalence classes.\newline
Let $R$ be such that $B_R(e)\cap S$ contains all elements of the 
equivalence classes with finitely many elements and at least one element of each 
equivalence class with infinitely many elements.\newline
For each vertex $y\in S$ we pick $t\in B_R(e)\cap S$, $t\sim y$ and
we consider the group $J$ generated by $\{g_y^{-1}g_t \}$.  Clearly $S\subset J(B_R(e)\cap S)$.
We claim that $Je$ is contained in a neighborhood of $S$. Indeed for each $g\in J$
we have that $gC_j$ and $C_j$ are at finite distance from each other for all $j$.
Let $A>0$ be such that for any $v\in \Gamma _H$ we have that if $d(v,S)>A$ then
$v$ lies in a deep component of $\Gamma _H -S$. One sees easily that such an $A$
exists by prop. 1, part 2.\newline
If $Je$ is not contained in any neighborhood of $S$ then there is a $g\in J$ such that
$gS\cap N_A(S)=\emptyset $. This follows again by lemma 2.2 and the remark at the end of the 
proof of lemma 2.2.
 Clearly then $S$ intersects a single component of
$\Gamma _H-gS$.
Therefore there is some $C_j$ such that all $gC_j$
except one
are contained in a single deep component of $ \Gamma _H-S$. This however  contradicts
the fact that $gC_j$ is at finite distance from $C_j$ for all $j$.\newline
We have therefore shown that $J$ is quasi-isometric to $S$. Clearly
$\Gamma _H/J$ has more than two ends. By the algebraic torus theorem of Dunwoody-Swenson
(\cite {D-S}) we have that $H$ splits over a group commensurable with $\Bbb Z ^n$.
In fact by proposition 3.1 of (\cite {D-S})
we have that $H$ splits over
a group commensurable with $J$. To see this note that if $T$ is an essential
track corresponding to $H$ as in lemma 2.3 of \cite {D-S} then from lemma 2.2
it follows that no translate of $T$ crosses $T$. 
$\blacksquare $
\enddemo
We return now to the proof of the theorem.
We proceed  inductively: We assume that we can write
$H$ as the fundamental group of a graph of groups $H=\pi _1(\Delta ,\B )$ where
 edge groups are commensurable with $J$ and $(\Delta ,\B )$ can not
be further refined. By this we mean that no vertex group $B_v$ of $\Delta $
admits a graph of groups decomposition with edge groups commensurable to $J $
and such that all edge groups $B_e $ of edges adjacent to $v$ are subgroups
of vertex groups of the graph of groups decomposition of $B_v$.\newline
By the accessibilty theorem of Bestvina-Feighn (\cite {B-F}) this procedure terminates.
Each vertex group of the graph of groups decomposition, say $H=\pi _1(\Delta ,\B )$,
that we obtain in this way
contains a finite index subgroup of $J$.
 Indeed if $|B_v:B_v\cap J|=\infty $ for some
vertex group then $B_v\cap J $ coarsely separates the Cayley graph of $B_v$, so by
applying the alebraic torus theorem once again we see that we can refine $(\Delta ,\B )$,
a contradiction.
We conclude therefore that all edge and vertex groups are commensurable
to $J$ which proves the theorem.
$\blacksquare $
\enddemo
\heading {\bf \S 3. Graphs of groups} \endheading
To deal with splittings over coarse $PD(n)$ groups in general rather than
$\Bbb Z ^n$ we need a theorem of Scott-Swarup (\cite {S-S}) generalizing proposition 3.1
of \cite {D-S}. We recall their notation
and results:
\proclaim {Definitions}
Two sets $P,Q$ are almost equal if their symmetric difference is finite.\newline
If $G$ acts on the right on a set $Z$ a subset $P$ is almost invariant
if $Pg$ is almost equal to $P$ for all $g\in G$.\newline
If $G$ is a finitely generated group and $J$ is a subgroup then
a subset $X$ of $G$ is $J$-almost invariant if $gX=X$ for all $g\in J$
and  $J\backslash X$ is an almost invariant subset of $J\backslash G$.\newline
$X$ is a non-trivial $J$-almost invariant subset of $G$ if, in addition,
$J\backslash X$ and $J\backslash (G-X)$ are both infinite.
\endproclaim
If $G$ splits over $J$ there is a $J$-almost invariant subset $X$ of $G$
associated to the splitting in a natural way. If $G=A\ast _J B$
let $X_A,X_B,X_J$ be complexes with 
$\pi _1(X_A)=A, \pi _1(X_B)=B, \pi _1(X_J)=J$ and let $X_G$ be
the complex obtained as usual by gluing $X_J\times I$ to $X_A\cup X_B$,
so that $\pi _1(X_G)=G$. Let $\tilde X_G$ be the universal covering
of $X_G$ and let $p:\tilde X_G\to X_G $ be the covering projection.
 We note that each connected component of $p^{-1}(X_J\times (0,1))$
separates $\tilde X$ in two sets. Each of these two sets gives
a $J$-almost invariant subset.\newline
If $X,Y$ are non-trivial $J$-invariant subsets of $G$ we say that
$Y$ {\it crosses } $X$ if all the four intersections
$X\cap Y,X^{\ast }\cap Y,X\cap Y^{\ast },X^{\ast }\cap  Y^{\ast }$
(where $ X^{\ast }=G-X, Y^{\ast }=G-Y$) project to infinite sets
in $J\backslash G$. If $X$ is a non-trivial $J$-invariant subset
we say that the intersection number  $i(J\backslash X, J\backslash X)$ is
$0$ if $gX$ does not cross $X$ for any $g\in G$. With this notation
Scott and Swarup (\cite {S-S}) show the following: 
\proclaim {Theorem  }
Let $G$ be a finitely generated group with a finitely generated subgroup
$J$ such that $\Gamma _G/J$ has more than one end. If there is a non-trivial
$J$ almost invariant set $X$ of $G$ such that $i(J\backslash X, J\backslash X)=0$
then $G$ has a splitting over some subgroup $J'$ commensurable with $J$.
\endproclaim
We note that from theorem 7.7 of \cite {K-K} it follows that if 
$f:X\to Y$ is a uniform embedding of a coarse $PD(n)$-space $X$
to a coarse $PD(n)$-space  $Y$ then $Y$ is coarsely contained in $f(X)$.
This implies that if $h:G\to H$ is a monomorphism between coarse $PD(n)$-groups then 
$|H:h(G)|<\infty $. \newline
\proclaim {Theorem 3.1 }
Let $G$ be a finitely generated group admitting a graph of groups decomposition
$G=\pi _1(\Gamma ,\A )$ such that all edge and vertex groups are coarse $PD(n)$ groups
of dimension $n$.
Suppose further that $(\Gamma,\A) $ is not a loop with all edge to vertex
maps isomorphisms and that it is not a graph of one edge with both edge to
vertex maps having as image an index 2 subgroup of the vertex group.
 If $H$ is quasi-isometric to $G$ 
then $H$ splits over a group that is quasi-isometric to 
 an edge group of $(\Gamma,\A) $.    
\endproclaim
\demo{Proof}
We  recall from \cite {S-W} the topological point of view on graphs of groups:
To each vertex $v\in \Gamma $ and each edge $e \in \Gamma $ 
 we associate  finite simplicial complexes $X_v,X_e$ such that
$\pi _1(X_v)=A_v, \pi _1(X_e)=A_e$. Let $I$ be the unit interval.
We construct a complex $X$ such that $\pi _1 (X)=G$
by gluing the complexes $X_v$ and $X_e\times I$ as follows:
Let $v$ be an endpoint of $e\in \Gamma $, say $v=\partial _0e$. 
Then there is a monomorphism $i_0:A_e\to A_v$. Let $f:X_e\to X_v$
be  a simplicial map such that $f_{\ast }=i_0$.
We identify then $(t,0)\in X_e\times I$ to $f(t)\in X_v$. \newline
Similarly we define an identification between $X_e\times \{1 \}$ and
$X_v$ if $v=\partial _1 e$.\newline
Doing all these identifications for the vertices $v\in \Gamma $ and the 
edges $e\in \Gamma $ we obtain a complex $X$ such that $\pi _1 (X)=G$.
We metrize the 1-skeleton of $\tilde X$ as usual by giving each edge length 1.
We note that with this metric $\tilde X$ is quasi-isometric to the
Cayley graph of $G$.\newline
We can obtain the Bass-Serre tree $T$ associated to $\pi _1(\Gamma ,\A )$ from $\tilde X$, by
collapsing each copy of $\tilde X_v\subset X $ to a vertex and each copy of
$\tilde X_e \times I$ to $I$. This collapsing gives a $G$-equivariant map 
$p:\tilde X\to T $.\newline
 We note now that if $v\in T$ is a vertex then $p^{-1}v$  separates $\tilde X $ into
more than 2 deep components. Moreover if $l$ is an infinite geodesic in $T$
$p^{-1}l$ is a coarse $PD(n+1)$-space (see theorem 11.13 of \cite {K-K}).
 As in lemma 2.2 we have
that if $h:\tilde X\to \tilde X $ is a quasi-isometry then for any vertex $v\in T$
there is a vertex $u\in T$ such that $h(p^{-1}(v))$ and $p^{-1}(u)$ are at finite distance
from each other.\newline
 
 Let $\Gamma _H$ be the Cayley graph of $H$.
Arguing as in lemma 2.3 we show that there is a subgroup $J$
of $H$ quasi-isometric to an edge group of  $(\Gamma ,\A )$ such that
$\Gamma _H/J$ has more than one end. As in lemma 2.3 we show
that there is a connected subset $S$ of $\Gamma _H$ such that
$ \Gamma _H-S$ has more than 1 deep components and $S/J$ is finite.
Without loss of generality we can assume that for any $v\in S$, $Jv\subset S$.
As in lemma 2.3 we see that we can find $J$ as above such
that ,in addition, $hC_i$ and $C_i$ are at finite distance from each other
for all deep components $C_i$ of $ \Gamma _H-S$. We fix now $v\in \Gamma _H$
and we identify $H$ with the orbit $Hv$. The set $X=J(C_i\cap Hv)$ is clearly
a non-trivial $J$ almost invariant subset of $H$. It is easy to verify
that $i(J\backslash X, J\backslash X)=0$ and  theorem 4.1 follows
from the result of Scott-Swarup quoted above.
$\blacksquare $.
\enddemo
\proclaim {Corollary 3.2  (\cite {F-M1,2}) }
Let $G$ be a solvable Baumslag-Solitar group.  If $H$
is a group quasi-isometric to $G$ then $H$ is commensurable
to a solvable Baumslag-Solitar group.
\endproclaim
\demo {Proof }
By theorem 3.1 $H$ splits over a 2-ended group. Since $H$ is amenable
and is not virtually abelian $H$ can be written as a graph of groups
with a single vertex and a single edge such that the edge group
is two ended and exactly one edge to vertex map is an isomorphism.
If $a$ is an element of the edge group generating an infinite normal
subgroup of the edge group and $t$ is the generator corresponding
to the edge then $tat^{-1}=a^k$ for some $k\in \Bbb Z$ and
the solvable Baumslag-Solitar subgroup of $H$ generated by $<t,a>$ is
a subgroup of finite index.
\enddemo
\proclaim {Corollary 3.3 }
Let $G$ be a finitely generated group admitting a graph of groups decomposition
$G=\pi _1(\Gamma ,\A )$ such that all edge and vertex groups are virtually $\Bbb Z ^n$.
Suppose further that $(\Gamma,\A) $ is not a loop with all edge to vertex
maps isomorphisms and that it is not a graph of one edge with both edge to
vertex maps having as image an index 2 subgroup of the vertex group.
 If $H$ is quasi-isometric to $G$ 
then $H=\pi _1(\Delta ,\Cal B)$ where all vertex and edge groups of $\Delta $
are virtually $\Bbb Z ^n$.
\endproclaim
\demo {Proof }
By theorem 3.1 $H$ splits over a subgroup quasi-isometric to $\Bbb Z ^n$
and hence virtually $\Bbb Z ^n$.
 We  apply the same argument as in the proof
of theorem 2.1 to conclude that  $H=\pi _1(\Delta ,\B )$ where all vertex and edge
groups of $\Delta $ are virtually  $\Bbb Z ^n$. $\blacksquare $

\enddemo
\proclaim {Definition }
A metric space $Y$ is dominated by a coarse $PD(n)$-space if
there is a  uniform embedding $f:Y\to X$ where $X$ is a coarse $PD(n)$-space.
We say that a finitely generated group $G$ is dominated by a coarse $PD(n)$-space if
$G$ equiped with the word metric is dominated by a coarse $PD(n)$ space.
\endproclaim
Some examples: A coarse $PD(n-k)$ group is dominated by a coarse $PD(n)$ space
where $k=0,...,n-1$. A free group is dominated by a coarse $PD(2)$ space.
\proclaim {Theorem 3.4 }
Let $G$ be a finitely generated group admitting a graph of groups decomposition
$G=\pi _1(\Gamma ,\A )$ such that all vertex groups are coarse $PD(n)$-groups
 and all edge groups are dominated by coarse $PD(n-1)$ spaces.
 If $H$ is quasi-isometric to $G$ 
then $H$ splits over some group quasi-isometric to an edge group
of $\Gamma $.
\endproclaim
\demo {Proof }
 We construct a complex $X$ such that $\pi _1(X)=G$ as in theorem 3.1. 
We have as in 3.1 that there is a map $p:\tilde X\to T$ where $T$
is the Bass-Serre tree of $(\Gamma ,\A )$. We note that $\tilde X$
is contained in a neighborhood of the `vertex spaces of $\tilde X$': 
 $p^{-1}(v)$ is a coarse $PD(n)$ space and $\tilde X$ is contained in the 1-neighborhood
of $\bigcup p^{-1}(v)$ (where $v$ ranges over vertices of $T$). Moreover
if $v,w$ are  vertices of $T$ adjacent to an edge $e$ then for $t>1$
 $N_t(p^{-1}(v))\cap N_t(p^{-1}(w))$
is a path connected subset of $\tilde X$ quasi-isometric to $p^{-1}e$ (here
we consider this subset equiped with its path metric). In fact $p^{-1}e$
and $N_t(p^{-1}(v))\cap N_t(p^{-1}(w))$ are contained in a finite neighborhood
of each other.
\newline
Let $f:\tilde X \to \Gamma _H$ be a quasi-isometry from $\tilde X$ to the Cayley
graph of $H$.
If $e$ is an edge of $T$ then $f(p^{-1}e)$ coarsely separates $\Gamma _H$.
Let $e$ be an edge of $T$ adjacent to the vertices $v,w$.
 Let $R_0$ be such that 
$f(p^{-1}e)\subset N_{R_0}(f(p^{-1}v))\cap N_{R_0}(f(p^{-1}w))$. 
 Note that
we can pick $R_0$ uniformly for all $e\in T$.
We distinguish now two cases:\newline
Case 1: $p^{-1}e$ is not quasi-isometric to a coarse $PD(n-1)$-space.\newline
Let $r$ be such that $f(p^{-1}v),f(p^{-1}w)$ are coarsely contained
in distinct components of $\Gamma _H- N_r(f(p^{-1}e))$. Let's call $F_1,F_2$ respectively these
2 components.
 Again
we can pick $r$ uniformly for all $e\in T$.
We set $S=N_r(f(p^{-1}e)), C_1=f(p^{-1}v), C_2=f(p^{-1}w) $.
Let $R_1'>2R_0$ be such that the following holds:
Let $u$ be a vertex of $T$ , $ C=hf(p^{-1}u)$ (where $h\in H$) and let $F$ be the component
 of $\Gamma _H-S$ coarsely containing $C$. 
Then if $x\in C-F$ we have $d(x,S)<R_1'$. Note that the existence of an $R_1'$ with this
property follows from prop. 1, part 2.\newline
Let $R_1>R_1'$ be such that the following holds:
If $x\in S$ then neither $B_x(R_1)\cap C_1$ nor $B_x(R_1)\cap C_2$ is contained
in the $R_1'$-neighborhood of $S$.\newline
Again we can pick $R_1',R_1$ uniformly for all $e,u\in T,h\in H$.\newline
Case 2: $p^{-1}e$ is quasi-isometric to a coarse $PD(n-1)$-space.\newline
We pick in this case too constants $r,R_1',R_1$ with similar properties:\newline
Let $r$ be such that $f(p^{-1}v),f(p^{-1}w)$ are not contained coarsely 
in the same component of $\Gamma _H- N_r(f(p^{-1}e))$. We pick $r$ so that
 $f(p^{-1}v), f(p^{-1}w)$ intersect each 2
 deep components of $\Gamma _H-N_r(f(p^{-1}e))$ along a half coarse $PD(n)$ space.
We set $S=N_r(f(p^{-1}e)), C_1=f(p^{-1}v), C_2=f(p^{-1}w) $.
Let $R_1>2R_0$ be such that the following hold:
If $u$ is a vertex of $T$, $ C=hf(p^{-1}u)$ (where $h\in H$) and $F_1,F_2$ are distinct
 components
 of $\Gamma _H-S$ such that $C\cap F_1,C\cap F_2$ are coarse half $PD(n)$-spaces then
if $x\in C-(F_1\cup F_2)$ then $d(x,S)<R_1'$. 
We suppose further that the following holds: If $F$ is any deep component of
$\Gamma _H-S$ that intersects $C_1$ (or $C_2$) along a half coarse $PD(n)$-space
 then $B_x(R_1)\cap F\cap C_1$ is not contained in the $R_1'$-neighborhood of $S$
 (and similarly for $C_2$.)
Again we can pick $r,R_1',R_1$ uniformly for all $e,u\in T,h\in H$.\newline
We note that there is an $R_2>R_1$ such that the following hold: \newline
a) for any $x\in S$  any two points in
$B_x(R_1)$ that lie in the same deep component can be joined by a path
lying in $B_x(R_2)-S$. \newline
b) for any $v\in B_x(R_1)$ $d(v,S)=d(v,S\cap B_x(R_2))$.
\newline
 We fix a vertex $v\in S$. For any vertex $x\in S$ we pick a $g_x\in H$
such that $g_xx=v$. We call the set of vertices in $B_v(R_2)\cap g_xS$
the {\it type } of $x$.
We will show the following:
\proclaim {Lemma 3.5} There is an $M>0$ such that if $x,y$ are of the same type then
$g_xS$ and $g_yS$ lie in the $M$-neighborhood of each other.
\endproclaim
\demo {Proof}
We note that two points in $B_v(R_1)$ lie in the same deep component     
of $\Gamma _H-g_xS$ if and only if they lie in the  same deep component
of $\Gamma _H-g_yS$. \newline
We distinguish two cases:\newline
Case 1: $S$ is not quasi-isometric to a coarse $PD(n-1)$-space.\newline
Let  $C_1=f(p^{-1}v), C_2=f(p^{-1}w) $ where $v,w$ are vertices adjacent to $e$.
$g_xC_1, g_xC_2$ are coarsely contained in distinct components, say $F_1,F_2$
of $\Gamma _H-g_xS$. Let $c_1\in B_v(R_1)\cap g_xC_1$ such that $d(c_1,g_xS)>R_1'$.
Then $c_1$ lies in $g_xF_1$. Since $x,y$ are of the same type $d(c_1,g_yS)>R_1'$
so $c_1$ lies in the  deep component of $\Gamma _H-g_yS$ that coarsely contains
$g_xC_1$. We pick similarly $c_2\in B_v(R_1)\cap g_xC_2$. Since $c_1,c_2$ are not contained
in the same component of $\Gamma _H-g_yS$ we have that $g_xC_1,g_xC_2$ are coarsely contained
in distinct deep components of $\Gamma _H-g_yS$.
 We claim that $N_{R_0}(g_xC_1)\cap N_{R_0}(g_xC_2)$ is contained
in the $R_1+R_0$ neighborhood of $g_yS$ so $g_xS$ is coarsely
contained in $g_yS$. Indeed let $a\in N_{R_0}(g_xC_1)\cap N_{R_0}(g_xC_2)$.
Let $a_1\in g_xC_1$, $a_2\in g_xC_2$ with $d(a,a_1)\leq R_0$, $d(a,a_2)\leq R_0$.
If $a_1$ (or $a_2$) is not contained in the deep component of $Gamma _H-g_yS$ that contains
$g_xC_1$ we have $d(a_1,g_yS)\leq R_1'$ so $d(a,g_yS)\leq R_0+R_1'$.
Otherwise we have that $a_1,a_2$ lie in distinct components of $Gamma _H-g_yS$
and there is a path joining them of length less or equal to $2R_0$. This path intersects
$g_yS$ so in this case $a$ is at distance less than $R_0$ from $g_yS$.
In the same way we see that $g_yS$ is coarsely
contained in $g_xS$. Clearly if $M=R_1+R+0$ 
we have that $g_xS$ and $g_yS$ lie in the $M$-neighborhood of each other.
We note in particular that $M$ does not depend on $x,y$.\newline

Case 2: $S$ is  quasi-isometric to a coarse $PD(n-1)$-space.\newline
Let  $C_1=f(p^{-1}v), C_2=f(p^{-1}w) $ where $v,w$ are vertices adjacent to $e$.
$g_xC_1, g_xC_2$ are not  contained coarsely in the same component 
of $\Gamma _H-g_xS$. Say $g_xC_1, g_xC_2$ intersect respectively
the (distinct) deep components $F_1,F_2$ of $\Gamma _H-g_xS$ along half coarse $PD(n)$ spaces.
Let $c_1\in g_xC_1\cap B_v(R_1)$ and $c_2\in g_xC_2\cap B_v(R_1)$
such that $d(c_1,g_xS)>R_1', d(c_2,g_xS)>R_1'$.
As in case 1 we have that $c_1,c_2$ lie in distinct components of $\Gamma _H-g_yS$
and so
$g_xC_1\cap F_1,g_xC_2\cap F_2 $ are
coarsely contained in distinct components of $\Gamma _H-g_yS$.
We have then that $N_{R_0}(g_xC_1)\cap N_{R_0}(g_xC_2)$ is contained
in the $R_1'+R_0$ neighborhood of $g_yS$ so $g_xS$ is coarsely
contained in $g_yS$. In the same way we see that $g_yS$ is coarsely
contained in $g_xS$. So if $M=R_1'+R_0$
we have that that $g_xS$ and $g_yS$ lie in the $M$-neighborhood of each other. $\blacksquare $ 
\enddemo
We consider now the group $J_1$ generated by all elements $\{g_x^{-1}g_y \}$ where
$x,y$ are of the same type.
We claim
that it follows from lemma 3.5 that
there is an $M_1>0$ such that for any $h\in J_1$, if $hS\cap S\ne \emptyset $
then $hS$ is contained in the $M_1$-neighborhood of $S$. \newline
Indeed we remark that there are finitely many types of vertices of $S$. Let
$x_1,...,x_n$ be vertices of $S$ representing these types and let $g_{x_1},...,g_{x_n}$
be the corresponding group elements. There is a $M_0>0$ such that for any $i,j\in \{1,..,n\}$ if
$g_{x_i}S $ and  $g_{x_j}S $ are contained in a finite neighborhood of each other then
they are contained in the $M_0$-neighborhood of each other.
Consider now an $h\in J_1$ such that $hS\cap S\ne \emptyset $. Since every generator 
of $J_1$ maps $S$ in a finite neighborhood of itself we have that $hS$ and $S$ are contained
in a finite neighborhood of each other.\newline
Similarly we see that for any $K>0$ there is a $M_K>0$ such that
for any $h\in J_1$ if $hS$ intersects the $K$-neighborhood of $S$ then
$hS$ is contained in the $M_K$-neighborhood of $S$.\newline
Let $a\in hS\cap S$. We have that $g_ahS$ and $g_aS$ are contained in a finite neighborhood
of each other. Therefore they are contained in an $M_0+2M$ neighborhood of each other.
Hence we can take $M_1=M_0+2M$.\newline
We fix $R>0$ such that $S\subset J_1(B_v(R)\cap S)$. We distinguish two cases:\newline
Case 1:There is an $M_2>0$ such that for any $h\in J_1$, $hS$ is contained
in the $M_2$-neighborhood of $S$. Then clearly $\Gamma _H/J_1$ has more
than 1 end. Let $S'=J_1S$. Then $\Gamma _H-S'$ has finitely
many deep components, say $F_1,...,F_n$ and $S'/J_1$ is finite.
Moreover $J_1$ acts on the set of deep components by permutations.
Therefore there is a finite index subgroup of $J_1$ ,say $J_1'$,
such that $F_1$ is a non-trivial $J_1'$-almost
invariant subset of $H$. Clearly $hF_1$ does not cross $F_1$ for
any $h\in H$.
It follows by \cite {S-S} that $H$ splits
over a subgroup commensurable with $J_1$. \newline
Case 2: Suppose now that no such $M_2>0$ exists. It follows
that there is an $h_0\in J_1$ such that $h_0S$ does not intersect the
$R_1'$-neighborhood of $S$. This implies that $\Gamma _H-(S\cup h_0S)$
has at least 3 deep components. We can now apply the same argument as in lemma 2.3
to conclude that $H$ splits. We do this in detail here:\newline
Let $A>M$ be such that $h_0S\subset N_A(S)$. We set $S'=N_A(S)$.
$\Gamma _H-S'$ has at least 3 deep components.\newline
 For each $x\in S$ we pick $h_x\in J_1$ such that $h_xx\in B_v(R)\cap S$. If 
$F_1,...,F_k$ are the deep components of $\Gamma _H -S$ and
$D_1,...,D_m$ the deep components of $\Gamma _H -S'$ we have
that for each $j$ there are $D_{j_1},...,D_{j_r}$ such that $h_xF_j$ and
$D_{j_1}\cup...\cup D_{j_r}$ are at finite distance from each other.\newline
We use this to define a map $f_x:\{ F_1,...,F_k \} \to \Cal P (\{ D_1,...,D_m \})$ 
where $f_x(F_j)=\{D_{i_1}, ..., D_{i_r} \}$ if and only if
$h_xF_j$ and $D_{i_1}\cup ...\cup D_{i_r}$ are at finite distance from each other. \newline
Let's write now $x\sim y$ for $x,y$ vertices of $S$ if $f_x=f_y$. Clearly $\sim $
is an equivalence relation with finitely many equivalence classes.\newline
Let $R_3$ be such that $B_{R_3}(v)\cap S$ contains all elements of the 
equivalence classes with finitely many elements and at least one element of each 
equivalence class with infinitely many elements.\newline
For each vertex $y\in S$ we pick $t\in B_R(v)\cap S$, $t\sim y$ and
we consider the group $J$ generated by $\{h_y^{-1}h_t \}$.  Clearly $S\subset J(B_{R_3}(v))$.
We note that for any $h\in J$ $hF_i$ and $F_i$ are in a finite neighborhood
of each other for all $i$.\newline
We claim that there is an $M_3>0$ such that for any $h\in J$, $hS$ is contained
in the $M_3 $-neighborhood of $S$.\newline
 Indeed, if for some $h\in J$, $hS$ does not intersect the $R_1'$-neighborhood of $S$
then the following hold: $hS$ lies in a single component of $\Gamma _H-S$, say $F_i$.
Moreover at least 2 deep components of $\Gamma _H-S$ are coarsely contained in $F_i$.
To see this recall that if $S$ is not a coarse $PD(n-1)$-space there are 2 deep components of
$\Gamma _H-S$, say $F_1,F_2$ containing, respectively, coarse $PD(n)$-spaces $C_1,C_2$,
such that for any $x\in C_1\, (C_2)$ if $d(x,S)\geq R_1'$ then $x$ lies in the deep
component of $\Gamma _H-S$ that coarsely contains $C_1\, (C_2)$. This implies that
$hC_1,hC_2$ are both contained in the same deep component of $\Gamma _H-S$, so $hF_1,hF_2$
are both contained in a finite neighborhood of a single deep component of $\Gamma _H-S$,
a contradiction.\newline
We argue similarly when $S$ is a coarse $PD(n-1)$-space. It follows that there is
a $M_3>0$ such that for any $h\in J$, $hS$ is contained
in the $M_3 $-neighborhood of $S$. 
Therefore the set $JS$ is at finite distance from $S$ , $JS/J$ is finite
and $\Gamma _H-JS$ has more than one deep component. Using the criterion of \cite {S-S} 
as before we conclude that $H$ splits over a group commensurable with $J$.
 $\blacksquare $
 
\enddemo
One can get a finer result when all edge groups are virtually $\Bbb Z ^{n-1}$:
\proclaim {Theorem 3.6 }
Let $G$ be a finitely generated group admitting a graph of groups decomposition
$G=\pi _1(\Gamma ,\A )$ such that all vertex groups are coarse $PD(n)$-groups
 and all edge groups are virtually $\Bbb Z ^{n-1}$.
 If $H$ is quasi-isometric to $G$ 
then $H=\pi _1(\Delta ,\B )$ where 
all edge
groups of $\Delta $ are virtually  $\Bbb Z ^{n-1}$ and
all vertex groups $H_v$ of $\Delta $ are either
coarse $PD(n)$-groups or virtually $\Bbb Z ^{n-1}$ or
 the pair $(H_v,\{ H_{e_i} \})$ is a coarse $PD(n)$-pair, where $H_{e_i}$
are the edge groups of the edges containing  $v$.

\endproclaim
\demo {Proof }
We construct a complex $X$ such that $\pi _1(X)=G$ as in theorem 3.4. 
We have as in 3.4 that there is a map $p:\tilde X\to T$ where $T$
is the Bass-Serre tree of $(\Gamma ,\A )$. 
  By theorem 3.4 $H$ splits over a group $J$ commensurable to $\Bbb Z ^{n-1}$. Moreover
$J$ lies at finite distance from $f(p^{-1}e)$ for some $e\in T$, where
$f:\tilde X \to \Gamma _H$ is a quasi-isometry. 
We proceed now inductively: We assume that we can write
$H$ as the fundamental group of a graph of groups $H=\pi _1(\Delta ,\B )$ where
 edge groups are commensurable with $\Bbb Z ^{n-1}$ lie at finite
distance from $f(p^{-1}e)$ for some $e\in T$
 and $(\Delta ,\B )$ can not
be further refined. By this we mean that no vertex group $B_v$ of $\Delta $
admits a graph of groups decomposition with edge groups of the same type
and such that all edge groups $B_e $ of edges adjacent to $v$ are subgroups
of vertex groups of the graph of groups decomposition of $B_v$.\newline
By the accessibilty theorem of Bestvina Freighn (\cite {B-F}) this procedure terminates.
We note now that if $H_v$ is a vertex group of $(\Delta ,\B )$ then
$f^{-1}(H_v)$ is coarsely contained in $p^{-1}u$ for some vertex $u\in T$. Indeed
if not we could refine  $(\Delta ,\B )$  further  as in proposition 3.3.\newline
If $H_v$ is not a coarse $PD(n-1)$-group and it is not quasi-isometric to
$p^{-1}u$ we have that $H_v$ acts by quasi-isometries (via $f$)
on the coarse $PD(n)$-space $p^{-1}u$ and it easily follows that 
  the pair $(H_v,\{ H_{e_i} \})$ is a coarse $PD(n)$-pair, where $H_{e_i}$
are the edge groups of the edges containing  $v$. $\blacksquare $
\enddemo
\heading {\bf \S 4. Questions } \endheading

It is reasonable to wonder whether theorems 3.1 and 3.4 can be subsumed under a theorem
posing no restriction on edge groups. More precisely we have the following:
\proclaim {Question 1 }
Let $G$ be a finitely generated group admitting a graph of groups decomposition
$G=\pi _1(\Gamma ,\A )$ such that all vertex groups are coarse $PD(n)$-groups.
 Is it true that
if $H$ is quasi-isometric to $G$ 
then $H$ splits over some group quasi-isometric to an edge group
of $\Gamma $?
\endproclaim
The main motivation of this paper was to generalize Stallings' theorem on groups
with infinitely many ends to splittings over groups that are not necessarily finite.
This has been achieved to a large extend for splittings over virtually cyclic groups
(see \cite { Bo}, \cite {P }). However generalizing Stallings' theorem for splittings
over any group poses serious difficulties. For example even in the virtually cyclic case
we have that surface groups do split over $\Bbb Z$ but triangle groups that are
quasi-isometric to them don't. In general there are examples of groups that are quasi-isometric
to groups that split but which do not split themselves and even do not virtually split 
(that is none
of their finite index subgroups splits). Moreover it is a hard problem to determine
whether a group virtually splits, it is not known even for 3-manifold hyperbolic groups.\newline
On the other hand looking closer at Stallings' theorem one realizes that it splits in
3 cases: Groups with 2 ends, virtually free groups and groups with more than 2 ends that
are not virtually free groups. Or to rephrase it in asymptotic topology terminology, the first two
cases correspond to groups of asymptotic dimension 1 which are coarsely separated by
subsets of  asymptotic dimension 0 (compact sets) and the third case corresponds
to groups of asymptotic dimension $\geq 2$ which are coarsely separated by
subsets of  asymptotic dimension 0 (compact sets). The first two cases can be considered
as `exceptional' as they belong to only 2 quasi-isometry (in fact commensurability) classes.
It seems that these two cases are the harder to generalize. On the other hand the
third case (the `codimension 2' case) might be easier to deal with. More precisely
we have the following:
\proclaim {Question 2 }
Let $G$ be a finitely generated group of asymptotic dimension $\geq n$.
Suppose that a uniformly embedded subset $S$ of asymptotic dimension $\leq n-2$
coarsely separates the Cayley graph of $G$. Is it true then that $G$ splits?
 \endproclaim
See \cite {G2} for a definition of asymptotic dimension.

\enddocument
\end
\bye